\documentclass[10pt,journal,letterpaper,compsoc]{IEEEtran}
\usepackage{amssymb,amsmath,url,epsfig,subfigure,graphicx,algorithm,algorithmic,times,flushend}

\begin{document}

\title{Preconditioned Krylov solvers for kernel regression}
\author{Balaji Vasan Srinivasan, Qi Hu, Nail A. Gumerov,
\IEEEcompsocitemizethanks{\IEEEcompsocthanksitem B.V. Srinivasan, Q. Hu, N.A. Gumerov and R. Duraiswami are with the Department of Computer Science, Institute for Advanced Computer Studies (UMIACS), University of Maryland, College Park, MD, USA\protect\\
E-mail: [balajiv,huqi,gumerov,ramani]@umiacs.umd.edu}
Raghu Murtugudde and
\IEEEcompsocitemizethanks{\IEEEcompsocthanksitem Raghu Murtugudde is with the Department of Atmospheric and Ocean Sciences, Earth Science Systems Interdisciplinary System (ESSIC), University of Maryland, College Park, MD, USA \protect\\
E-mail: ragu@essic.umd.edu}
Ramani Duraiswami

\thanks{This work is a part of Chesapeake Bay Forecast System and we gratefully acknowledge National Ocean and Atmospheric Administration (NOAA) Award NA$06$NES$4280016$ for funding the project. We also acknowledge NSF award 0403313 and NVIDIA support for the Chimera cluster at the CUDA Center of Excellence at UMIACS.}}

\markboth{IEEE Transactions on Pattern Analysis and Machine Intelligence}%
{Srinivasan \MakeLowercase{\textit{et al.}}: Preconditioned Krylov solvers for kernel regression}

\IEEEcompsoctitleabstractindextext{
\begin{abstract}
A primary computational problem in kernel regression is solution of a dense linear system with the $N\times N$ kernel matrix. Because a direct solution has an O($N^3$) cost, iterative Krylov methods are often used with fast matrix-vector products. For poorly conditioned problems, convergence of the iteration is slow and preconditioning becomes necessary. We investigate preconditioning from the viewpoint of scalability and efficiency. The problems that conventional preconditioners face when applied to kernel methods are demonstrated. A \emph{novel flexible preconditioner }that not only improves convergence but also allows utilization of fast kernel matrix-vector products is introduced. The performance of this preconditioner is first illustrated on synthetic data, and subsequently on a suite of test problems in kernel regression and geostatistical kriging.
\end{abstract}
}

\maketitle
\IEEEdisplaynotcompsoctitleabstractindextext
\IEEEpeerreviewmaketitle

\section{Introduction}
The basic computations in kernel regression methods involve a number of linear algebra operations on matrices of kernel functions ($\bf \hat{K}$), which take as arguments the training and/or the testing data. A kernel function $k(x_i,x_j)$ generalizes the notion of the similarity between data points. Given $\mathbf{X}=\{x_1,x_2,\ldots,x_N\}, x_i\in R^d$, the kernel matrix entries are given by,
\begin{equation}
\mathbf{\hat{K}} = \left(
  \begin{array}{ccc}
    k(x_1,x_1) & \ldots & k(x_1,x_N)\\
    \vdots & \ddots & \vdots\\
    k(x_N,x_1) & \ldots & k(x_N,x_N)\\
  \end{array}
\right).\label{eq:kernelMatrix}
\end{equation}
The kernel function $k$ is chosen to reflect prior information, in the absence of which, the Gaussian kernel $\Phi$ is widely used,
\begin{equation}
 \Phi(x_i,x_j)=\exp\left(-\frac{\|x_i-x_j\|^2}{2\sigma^2}\right)\label{eq:Gaussian}.
\end{equation}
We use this kernel, though the methods discussed apply to other kernels as well, as is illustrated in experiments. The kernel matrix is usually regularized,
\begin{equation}\mathbf{K}=\mathbf{\hat{K}}+\gamma\mathbf{I};\label{eq:RegKernel}\end{equation}
with $\gamma$ chosen appropriately according to the problem.

Kernel regression appears in many variations: e.g. ridge regression \cite{BishopML}, Gaussian process regression \cite{GPML_Rasmussen} and geostatistical kriging \cite{AppliedGeostatistics}. The key computation in all these variants is the solution of a linear system with $\mathbf{K}$.

Direct solution for a dense kernel matrix system has a time complexity O($N^3$) and a memory complexity O($N^2$), which prevents its use with large datasets. Iterative Krylov methods \cite{SaadIterativeMethods} address this partially by reducing the time complexity to O($kN^2$), $k$ being the number of iterations \cite{GPML_Mackay,NandoKrylov}. The dominant cost per Krylov iteration is a kernel matrix-vector product (MVP), whose structure has been utilized to reduce the O($N^2$) space and time complexity further. The space requirement is reduced to O($N$) by casting the MVP as a weighted kernel summation and computing $k(x_i,x_j)$ \emph{on-the-fly} when evaluating the sum. Further, by using efficient kernel MVP \cite{IFGT_Yang,IFGT_Raykar,DualTree,FIGTREE,GPUML}, the cost of the MVP in each Krylov iteration can be reduced to O($N\log N$) or O($N^2/p$), $p$ being the number of processors. \emph{In these fast kernel MVP, there is usually a trade-off between accuracy and speed, and usually a MVP of reduced accuracy can be obtained faster.} This is either explicit (e.g. single precision SSE or graphical processors \cite{GPUML}) or algorithmic (IFGT \cite{IFGT_Yang,IFGT_Raykar}, dual-tree \cite{DualTree}, Figtree \cite{FIGTREE}). However the convergence rate suffers as the problem size increases since the matrix condition number usually increases with data. To speedup iterative methods in these cases, apart from using fast MVP, we need to reduce the number of iterations.

The convergence of the Krylov methods is determined by the matrix condition number $\kappa$ ($\kappa\ge1$),
\begin{equation}
\kappa=\frac{\lambda_{\max}}{\lambda_{\min}}, 1\le\kappa<\infty.
\end{equation}
where $\lambda_{\max}$ and $\lambda_{\min}$ are the largest  and smallest eigenvalues of $\mathbf{K}$ respectively. For smaller $\kappa$, the convergence is faster. For larger $\kappa$, there is a significant decrease in the convergence rate, necessitating a ``preconditioner'' \cite{SaadIterativeMethods} to improve the conditioning. Preconditioning has been suggested for kernel methods \cite{NandoKrylov,MurrayKrylov}, but to our knowledge, there has been \emph{no previous work to design a preconditioner }for such matrices.

To be effective, the preconditioner matrix construction cost should be small, and be able to take advantage of fast MVPs. We propose a novel preconditioner that improves convergence and has the added benefit that it utilizes the fast MVPs available for the kernel matrix.

The paper is organized as follows. We discuss kernel regression and its variants that we seek to use in Sec. \ref{sec:kernelRegression}. We introduce Krylov methods and their convergence properties in Sec. \ref{sec:Krylov} and survey different preconditioning techniques in Sec. \ref{sec:precond}. The new preconditioner is introduced and its parameters and convergence are studied in Sec. \ref{sec:ProposedPreconditioner}. Finally we test its performance on synthetic and standard datasets in Sec. \ref{sec:Experiments}.

\section{Kernel regression\label{sec:kernelRegression}}
We are particularly interested in Gaussian process regression and geostatistical kriging.

\subsection{Gaussian process regression (GPR) \cite{GPML_Rasmussen} \label{sec:GPR}}
GPR is a probabilistic kernel regression approach which uses the prior that the regression function is sampled from a Gaussian process. Given $D=\{x_i,y_i\}_{i=1}^N$, where $x_i$ is the input and $y_i$ is the corresponding output, the function model is assumed to be $y=f(x)+\epsilon$, where $\epsilon$ is a Gaussian noise process with zero mean and variance $\gamma$. Rasmussen et al. \cite{GPML_Rasmussen} use the Gaussian process prior with a zero mean function and a covariance function defined by a kernel $\mathbf{\hat{K}}(x,x_\ast)$, which is the covariance between $x$ and $x_\ast$, i.e. $f(x)\sim GP(0,\mathbf{\hat{K}}(x,x_\ast))$. With this prior, the posterior of the output $f(x_\ast)$ is also Gaussian with mean $m$ and variance $\Sigma$:
\begin{equation}m = k(x_\ast)^{T}(\mathbf{\hat{K}}+\gamma I)^{-1}y, \label{eq:GPR_mean}\end{equation}
\begin{equation}\Sigma = \mathbf{\hat{K}}(x_\ast,x_\ast)-k(x_{\ast})^{T}(\mathbf{\hat{K}}+\gamma I)^{-1}k(x_{\ast})\end{equation}
where $x\ast$ is the input at which prediction is required and $k(x_\ast)=[\mathbf{\hat{K}}(x_1,x_\ast),\mathbf{\hat{K}}(x_2,x_\ast)\ldots,\mathbf{\hat{K}}(x_N,x_\ast)]^T$. Here ``inverses'' imply solution of the corresponding linear system. Hyper-parameters (eg. $\sigma$ in Eq. \ref{eq:Gaussian}) are estimated via maximum likelihood techniques \cite{GPML_Rasmussen}. Note that the noise variance $\gamma$ results in the regularization of the kernel matrix, hence is similar in its role to the one in Eq. \ref{eq:RegKernel}.

\subsection{Kriging }
Kriging \cite{AppliedGeostatistics} is a group of geostatistical techniques to interpolate the value of a random field at an unobserved location from observations of its value at nearby locations. It was first used with mining and has since been applied in several scientific disciplines including atmospheric science, environmental monitoring and soil management.

There are several versions of kriging; the commonly used \emph{simple kriging} results in a formulation similar to Gaussian process regression \cite{AppliedGeostatistics}. Given geostatistical values $y_i$s recorded at locations $x_i$s, the interpolation at a new point $x_\ast$ is given by,
\begin{equation}y_\ast=k(x_\ast) (\mathbf{\hat{K}}+\gamma \mathbf{I})^{-1}y,\label{eq:Kriging}\end{equation}
where $k(x_\ast)$ is similar to the posterior mean in Eq. \ref{eq:GPR_mean}.

\section{Krylov methods \label{sec:Krylov}}
Krylov methods are formulated as a ``cost-minimization'' problem over a set of basis vectors (the Krylov basis) created via matrix vector products of the matrix under consideration. A detailed discussion and analysis can be found in \cite{SaadIterativeMethods,SaadEigen}; we provide a brief overview here.

For solving $\mathbf{K}x=b$. Krylov methods begin with an initial guess $x^{(0)}$ and minimize the residual $r^{(k)}=b-\mathbf{K}x^{(k)}$ in some norm, by moving the iterates along directions in the Krylov subspace $\mathcal{K}_k=span(r_0,\mathbf{K}r_0,\ldots,\mathbf{K}^{k-1}r_0)$. The directions are augmented over each Krylov iteration, a significant difference from simpler iterative approaches like Gauss-Siedel where the next iterate depends only on the previous one.

At the $k^{th}$ iteration, an orthogonal matrix $V^{(k)}=\left[v_1,v_2,\ldots,v_k\right]$ is generated such that columns of $V^{(k)}$ span the Krylov subspace $\mathcal{K}_k$ \cite{SaadIterativeMethods},
\begin{equation}\mathbf{K}V^{(k)}=V^{(k+1)}\mathbf{\bar{H}}^{(k)},\label{eq:Arnoldi}\end{equation}
where $\mathbf{\bar{H}}^{(k)}$ is an augmented Hessenberg matrix,
 \begin{eqnarray}\mathbf{\bar{H}}^{(k)}=\left(\begin{array}{ccccc}
                                      h_{1,1} & h_{1,2} & h_{1,3} & \ldots & h_{1,k} \\
                                      h_{2,1} & h_{2,2} & h_{2,3} & \ldots & h_{2,k} \\
                                      \vdots  & \vdots & \vdots & \vdots & \vdots \\
                                      0       & \ldots & 0 & h_{k,k-1} & h_{k,k} \\
                                      0       & \ldots & 0 & 0         & h_{k+1,k}
                                    \end{array}\right),\nonumber\end{eqnarray}
where $h_{i,j}=(v_j^T\mathbf{K}v_i)$. The next iterate $x^{(k)}$ is then given by,
\begin{equation}x^{(k)}=V^{(k)}\hat{y},\label{eq:ArnoldiIterate}\end{equation}
where $\hat{y}$ is obtained by solving the least squares problem, $\min_{\hat{y}} \|\mathbf{\bar{H}}^{(k)}\hat{y}-\beta e_1\|; e_1=\left[1,0,\ldots,0\right]^T$. This is the \emph{Arnoldi} iteration for system solution \cite{SaadIterativeMethods}.

The \emph{conjugate gradient (CG)} method is the most widely used Krylov method with symmetric matrices. For symmetric $\mathbf{K}$, $\mathbf{\bar{H}}^{(k)}$ in Eq. \ref{eq:Arnoldi} is tridiagonal making CG particularly efficient. The \emph{generalized minimum residual (GMRES)} is usually used for non-symmetric problems; GMRES minimizes the residuals $r^{(k)}$ in the $2-$norm. CG minimizes the $\mathbf{K}$-norm of the residual and utilizes the conjugacy in the resulting formulation, which results in not requiring to store the Krylov basis vectors. CG, therefore, is more efficient (lower cost per iteration) than GMRES. Kernel matrices are symmetric and satisfy the M\'ercer conditions $a^T\mathbf{K}a>0$, for any $a$; and hence $\mathbf{K}$ is positive definite. Therefore, when preconditioning is not used, CG has been the preferred choice \cite{GPML_Mackay}; however, GMRES has also been used \cite{NandoKrylov}.

The convergence rate is given by the ratio of the error ($e_k$) at $k^{th}$ iteration to the initial error ($e_0$) in some norm. For example, the ratio for CG \cite{SaadIterativeMethods} is,
 \begin{equation}\frac{\|e_k\|_\mathbf{K}}{\|e_0\|_\mathbf{K}}\le2\left(\frac{\sqrt{\kappa}-1}{\sqrt{\kappa}+1}\right)^k.\end{equation}
 A similar expression may be derived for GMRES \cite{SaadIterativeMethods}.

\subsection{Fast matrix-vector products: }
The key computation in each Krylov step is the MVP $\mathbf{Kq}$ or $\sum_{i=1}^{N}q_{i}k(x_{i},x_{j})$ for some vector $q$. Existing approaches to accelerate the MVP either approximate it \cite{IFGT_Yang,IFGT_Raykar,DualTree,FIGTREE} and/or parallelize it \cite{GPUML}; and have their pros and cons. In this paper, we present results with GPUML \cite{GPUML}, an open-source package that parallelizes kernel summation on graphical processors (GPUs) though we also tried with FIGTREE \cite{FIGTREE}. GPUML is easily extendable to generic kernels and works well with reasonable input data dimensions (up to $100$).

\subsection{Need for preconditioning: }
\begin{figure}[bth]
\centering
\includegraphics[width=\linewidth]{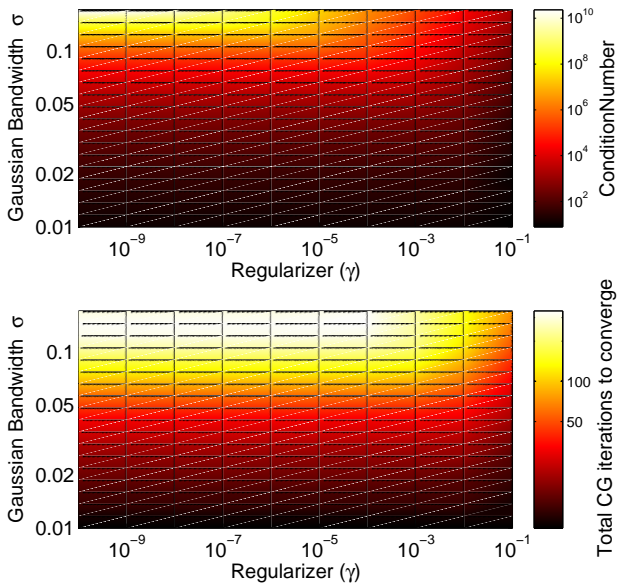}
\caption{\emph{Effect of kernel hyper-parameters on the matrix conditioning and CG iterations}\label{fig:GaussianHyperparameters}}
\end{figure}

The condition number $\kappa$ of kernel matrices depends on the data point distribution and the kernel hyper-parameters. For the Gaussian (Eq. \ref{eq:Gaussian}), the hyper-parameters are the bandwidth $\sigma$ and the regularizer $\gamma$. While $x_i$'s are given, the hyper-parameters are generally evaluated using ML. Fig. \ref{fig:GaussianHyperparameters} shows the $\kappa$ and number of CG iterations to converge for a kernel matrix constructed from data points chosen uniformly at random from inside a unit cube. We observe the following: there is a direct correspondence between $\kappa$ and number of CG iterations. For larger regularizer ($\gamma$) and smaller bandwidths ($\sigma$), the convergence is much better. The data point distribution influences the conditioning as well. It is however not possible to hand select these parameters for each problem. It is necessary to ``precondition'' \cite{SaadIterativeMethods} the system to be yield better convergence irrespective of the underlying matrix.

\section{Preconditioning techniques\label{sec:precond}}
A left preconditioner ($\mathbf{M}^{-1}$) operates on ${\bf K}x=b$ as,
\begin{equation}\mathbf{M}^{-1}\mathbf{K}x=\mathbf{M^{-1}}b; \label{eq:LeftPreconditioner}\end{equation}
and a right preconditioner operates as,
\begin{equation} \mathbf{KM}^{-1}y=b, \quad y=\mathbf{M}x.\label{eq:RightPreconditioner}\end{equation}
The Preconditioner $\mathbf{M}^{-1}$ should be chosen so that the $\mathbf{M}^{-1}\mathbf{K}$ or $\mathbf{K}\mathbf{M}^{-1}$ have a low $\kappa$. An ideal preconditioner ($\mathbf{M}^{-1}$) should well approximate $\mathbf{K}^{-1}$, but be easy to compute.

\subsection{Conventional preconditioners}
Standard preconditioners used in the literature were developed for sparse matrices that arise in the solution of differential equations, and include Jacobi and Symmetric Successive Over-Relaxation (SSOR). For general sparse matrices, incomplete LU or Cholesky preconditioners are often used. The triangular factors $\mathbf{L}$ and $\mathbf{U}$  for a sparse matrix may not be sparse, but {\em incomplete LU} factorizations leads to sparse $\mathbf{L}$ and $\mathbf{U}$ matrices by setting the coefficients leading to zero entries of the sparse matrix to zero. For a dense matrix, elements are sparsified using a cut-off threshold.

Preconditioners to radial basis function interpolation are a closely related problem. Fast preconditioners have been proposed \cite{RBF_Beatson, RBF_Faul, RBF_Gumerov}, however, these approaches are limited to low data dimensions ($\sim3$ dimensions for $X$).

\subsection{Flexible preconditioners}
As seen from Eqs. (\ref{eq:LeftPreconditioner}) and (\ref{eq:RightPreconditioner}), a left preconditioner modifies the right-hand side $b$ in the problem whereas the right preconditioner leaves it as is. This property of right preconditioners can be exploited to create ``flexible'' preconditioning techniques where a different preconditioner can be used in each Krylov step \cite{FGMRES,FlexibleKrylov,FlexibleCG}, since the preconditioner only appears implicitly. Flexible preconditioning can be used with both CG \cite{FlexibleCG} and GMRES \cite{FGMRES}.

Although many papers have shown the convergence of flexible preconditioners under exact arithmetic, it is very hard to estimate convergence rate or the number of outer iterations accurately under inexact arithmetic since the underlying subspaces, $x_0+span\{\mathbf{M_1}^{-1}v_1,\mathbf{M_2}^{-1}v_2, \ldots, \mathbf{M_k}^{-1}v_k\}$ are no longer a standard Krylov subspace. This affects CG since conjugacy is essential and cannot be guaranteed. Notay \cite{FlexibleCG} proposes $2$ modifications to a preconditioned flexible CG. The iterates should be ``reorthogonalized'' at each step to maintain conjugacy; and the preconditioner system should be solved with high accuracy. Flexible preconditioners are however easily used with GMRES. This fact will be observed in results below, where a poorer performance is observed for flexible CG relative to flexible GMRES.

The algorithmic details of flexible GMRES are shown in Algorithm \ref{alg:fgmres}, and the corresponding unpreconditioned version is obtained by replacing the $M$s with identity matrices. A similar extension is available for CG as well. The iterations are stopped when $\epsilon=\frac{b-\mathbf{K}x_i}{N}$ drops below a certain tolerance.

\begin{algorithm}[bth]
\caption{Flexible GMRES \cite{FGMRES}} \label{alg:fgmres}
\begin{algorithmic}[1]
\STATE $r_0=\left(b-\mathbf{K}x_0\right)$, $\beta=\|r_0\|_2$ and $v_1=r_0/\beta$
\STATE Define the $m+1 \times m$ matrix, $\bar{H}_m=\{h_{i,j}\}_{1\leq i\leq j+1;1\leq j\leq m}$
\FOR{$j=0$ to $iter$}
    \STATE Solve $\mathbf{M_j}z_j=v_j$ (\emph{\textbf{inner preconditioner}})
    \STATE $w=\mathbf{K}z_j$ (matrix-vector product)
    \FOR {$i=0$ to $j$}
        \STATE $h_{i,j} = (w,v_i)$, $w = w - h_{i,j}v_i$
    \ENDFOR
    \STATE $h_{j+1,j}=\|w\|_2$, $v_{j+1}=w/h_{j+1,j}$
\ENDFOR
\STATE $\mathbf{Z}_{(iter)}=\left[z_1,\ldots,z_{iter}\right]$,
\STATE $y_{iter}=\arg\min_y \|\beta e_1 - \bar{H}^{iter}y\|_2 $, $x_{iter} = x_0 + \mathbf{Z}_{iter}y_{iter}$
\STATE IF satisfied STOP, else $x_0=x_{iter}$ and GOTO $1$
\end{algorithmic}
\end{algorithm}

\subsection{Krylov method as a flexible preconditioner:} In Algorithm \ref{alg:fgmres}, all that is needed to prescribe the right preconditioner is a black-box routine which returns the solution to a linear system with the preconditioner matrix $\mathbf{M}$. Thus, instead of explicitly specifying $\mathbf{M}^{-1}$, it is possible to specify it implicitly by solving a linear system with $\mathbf{M}$ using another Krylov method such as CG. \emph{However, because this iteration does not converge exactly the same way each time it is applied, this is equivalent in exact arithmetic to using a different $\mathbf{M}$ for each iteration \cite{FlexibleKrylov}.} We refer to the preconditioner, operating with matrix $\mathbf{M}$ as ``inner Krylov'' and to the main solver with matrix $\mathbf{KM}^{-1}$ as ``outer Krylov''.

\section{Preconditioning kernel matrices \label{sec:ProposedPreconditioner}}
Conventional preconditioners require construction of the complete preconditioner matrix $\mathbf{M}$ initially, followed by expensive matrix decompositions. Thus they have a computational cost of O($N^3$) and a memory requirement of at least $O(N^2)$. Additionally, the preconditioner evaluations will require a O($N^2$) ``unstructured'' matrix-vector product, which does not have any standard acceleration technique and is harder to parallelize. This limits their application to very large datasets and will ruin any advantage gained by the use of fast matrix-vector products (as will be seen later in Sec. \ref{sec:TestOfConvergence}).

This leads us to propose a key requirement for any preconditioning approach for a kernel matrix: \emph{\textbf{the preconditioner should operate with an asymptotic time complexity and memory requirement that are at least the same as the fast matrix vector product.}} One of the main contributions of the paper is a particularly simple construction of a right preconditioner, which also has a fast matrix vector product.

We propose to use a \emph{regularized kernel matrix $\mathbf{K}$} as a right preconditioner,
\begin{equation}\mathbf{M}=\mathbf{K}+\delta \mathbf{I}.\end{equation}
Regularization is a central theme in statistics and machine learning \cite{VapnikStat}, and is not a new concept for kernel machines, e.g. ridge regression, where the kernel matrix ($\mathbf{\hat{K}}$) is regularized as $\mathbf{\hat{K}}+\gamma\mathbf{I}$. However, the $\gamma$ is chosen by statistical techniques, and hence cannot be controlled.

Our use of this old trick of regularization is in a new context -- in the preconditioner matrix $\mathbf{M}$. The simple prescription achieves the following properties:
\begin{itemize}
\item improves condition number of matrix ${\mathbf M}$, leading to faster convergence of inner iterations
\item improves conditioning of outer matrix $\mathbf{KM}^{-1}.$
\end{itemize}
To translate this idea in to a useful preconditioner, we need a prescription for selecting the regularization parameter $\delta$ and specifying the accuracy $\epsilon$ to which the inner system needs to be solved. Because CG is more efficient for unpreconditioned symmetric systems, we use it to solve the inner system.

\subsection{Preconditioner acceleration}
A preconditioner improves convergence at the expense of an increased cost per iteration. This is because there is a cost associated with the preconditioner construction (amortized over all iterations) and cost of the inner iteration. To be useful, the total time taken by the preconditioned approach should be smaller.

The key advantages of the proposed preconditioner is that, because $\mathbf{M}$ is derived from $\mathbf{K}$, given $\mathbf{X}=\{x_1,x_2,\ldots,x_N\}, x_i\in R^d$ it is not necessary to explicitly construct the preconditioner $\mathbf{M}^{-1}$. Further, the key computation in the inner Krylov iteration is a MVP, $\mathbf{M}x$. This can be accelerated using the same fast algorithm as for $\mathbf{K}$. Further, the preconditioner system only needs to be solved approximately (with a low residual CG tolerance and with a lower accuracy MVP). In our experiments we use low-accuracy fast matrix vector products for the inner iterations (single precision on the GPU). For the outer iterations, the products are performed in double-precision.

\subsection{Preconditioner parameters}
The preconditioner regularizer $\delta$ must be chosen on the one hand to converge quickly, while on the other hand not causing it to deviate too much from the inverse of ${\bf K}$. The convergence of CG for a kernel matrix for different $\delta$'s is shown in Fig. \ref{fig:GaussianHyperparameters}. It can be seen that for large enough $\delta$, CG converges rapidly. The CG can also be forced to have an early termination by setting a low solution accuracy ($\epsilon$).

\subsection{Effect of regularization parameter ($\delta$):} In flexible Krylov methods, the outer GMRES iteration solves $\mathbf{KM}^{-1}y=b$, and the inner CG solves $\mathbf{M}x=y$. For small $\delta$, $\mathbf{M}$ is closer to $\mathbf{K}$. Therefore, the outer iteration is better conditioned; however, when $\mathbf{K}$ is ill-conditioned, $\mathbf{M}$ will also be somewhat ill-conditioned, thus slowing the inner iterations.

\begin{figure}[p]
\centering
\subfigure[Effect of regularizer $\delta$ on flexible CG]{
\includegraphics[width=\linewidth]{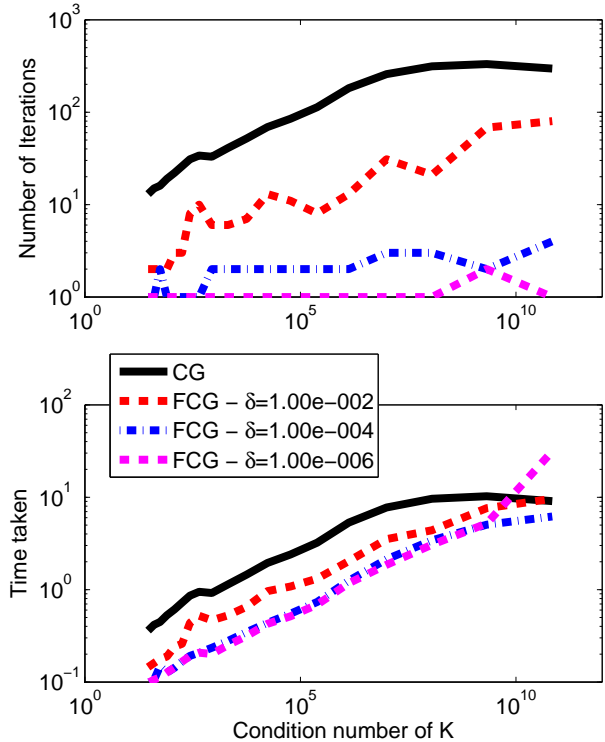}
\label{fig:fcg_sigma}}
\subfigure[Effect of regularizer $\delta$ on flexible GMRES]{
\includegraphics[width=\linewidth]{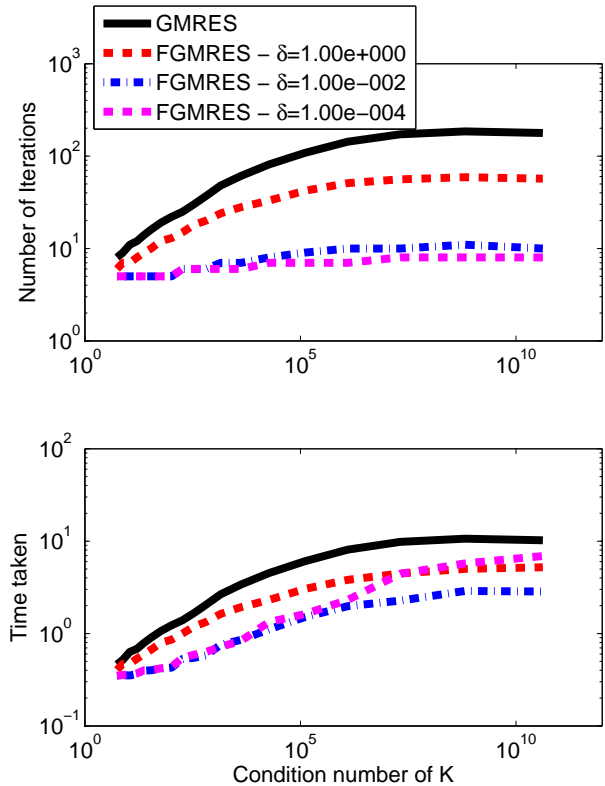}
\label{fig:fgmres_sigma}}
\caption{\emph{Effect of regularizer $\delta$ on the convergence for FCG and FGMRES for different conditioning of $K$. The condition number is adjusted by increasing the Gaussian bandwidth $\sigma$ for $\mathbf{K}$.}}
\end{figure}

To demonstrate this, we generated data as before by taking $2000$ random samples in a unit cube and generated a matrix for the Gaussian kernel. We tested the convergence with this preconditioner for various regularizer values (Figs. \ref{fig:fcg_sigma} and \ref{fig:fgmres_sigma}). For smaller $\delta$, the convergence of the outer iterations is faster, but the cost per iteration increases due to slow convergence of the inner iterations. Large regularization results in a poor preconditioner $\mathbf{M}$. An intermediate value of the regularizer is therefore optimal. This is observed for both flexible CG (FCG) and flexible GMRES (FGMRES). However, because of its formulation, the optimal FCG regularizer $\delta$ is an order of magnitude lower than that for FGMRES.

\emph{The choice of a regularizer involves a trade-off between the preconditioner's accurate representation of the kernel matrix and its desired conditioning.}

\subsection{Effect of CG tolerance ($\epsilon$):} 
\begin{figure}[p]
\centering
\subfigure[Effect of inner CG tolerance $\epsilon$ on flexible CG]{
\includegraphics[width=\linewidth]{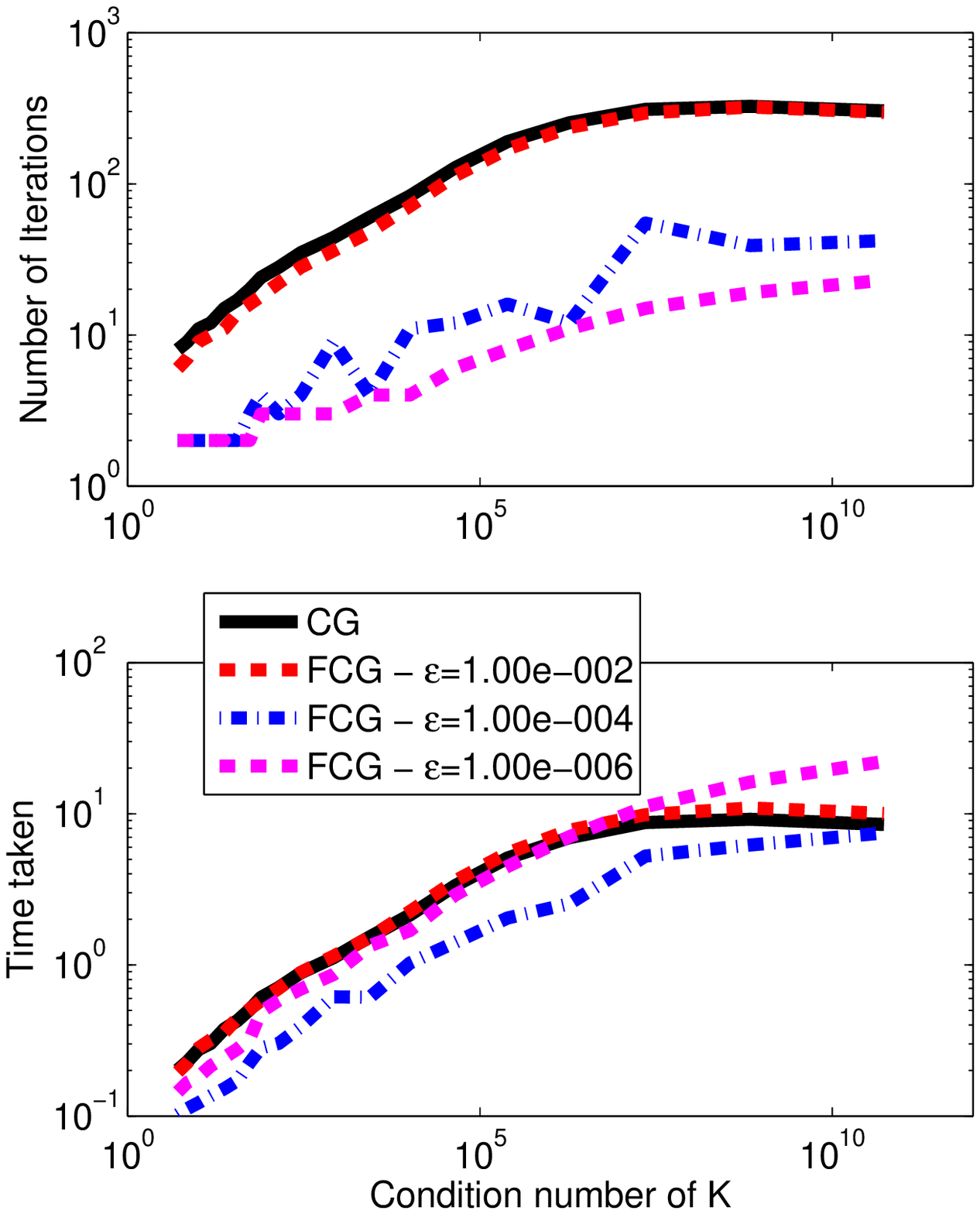}
\label{fig:fcg_tol}}
\subfigure[Effect of inner CG tolerance $\epsilon$ on flexible GMRES]{
\includegraphics[width=\linewidth]{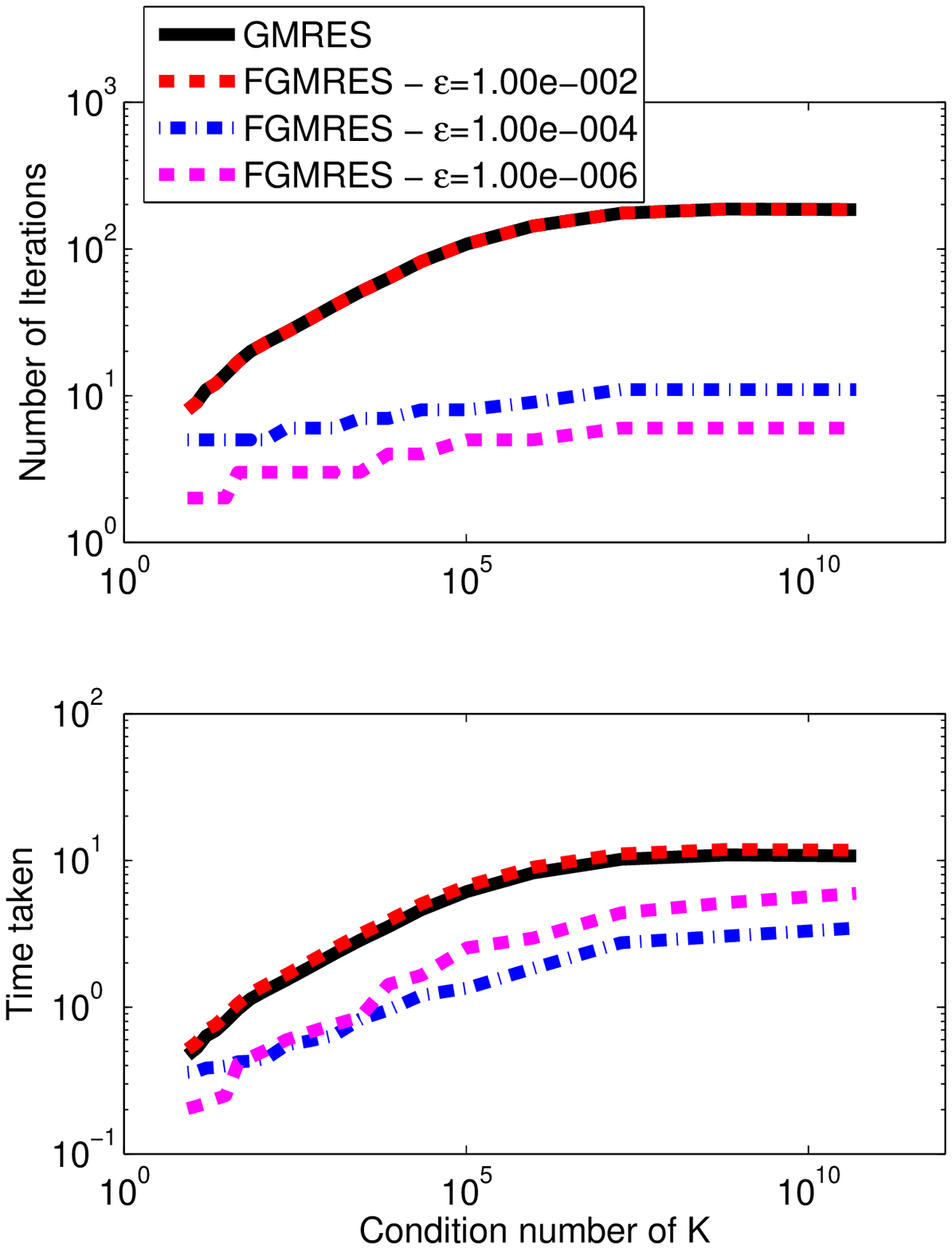}
\label{fig:fgmres_tol}}
\caption{\emph{Effect of CG tolerance $\epsilon$ on the convergence for FCG and FGMRES for different conditioning of $K$. The condition number is adjusted by increasing the Gaussian bandwidth $\sigma$ for $\mathbf{K}$.}}
\end{figure}

We tested the performance of the preconditioner for various tolerances in the inner iterations (Figs. \ref{fig:fcg_tol} and \ref{fig:fgmres_tol}). There is a consistent improvement in the outer convergence for more precise convergence settings of the inner solver. However, the cost of inner iterations increases. Therefore, an optimal intermediate value of $\epsilon$ works best for both FCG and FGMRES.

\emph{The choice of tolerance for CG iterations is a trade-off between the required solution accuracy of the   preconditioner system (and hence the convergence of the outer iterations) and the related computational cost. }

\subsection{Test of convergence\label{sec:TestOfConvergence}}
We compared the performance of FCG and FGMRES against ILU preconditioned CG and GMRES and the unpreconditioned CG and GMRES.

\begin{figure}[bth]
\centering
\includegraphics[width=\linewidth]{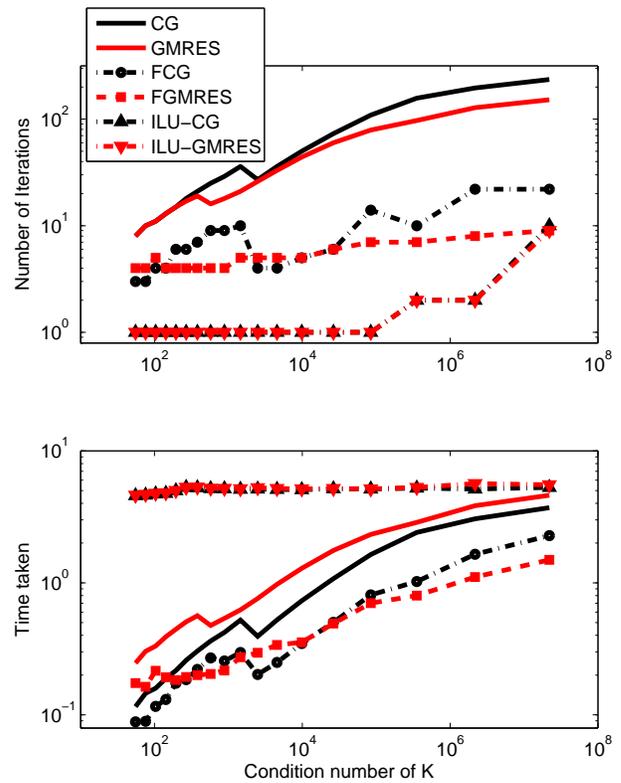}
\caption{\emph{Performance of our preconditioner with CG and GMRES against ILU-preconditioned and unpreconditioned versions for different conditioning of $K$. The condition number is adjusted by increasing the Gaussian bandwidth $\sigma$ for $\mathbf{K}$.\label{fig:fgmres}}}
\end{figure}

We set the preconditioner $\delta$ and tolerance $\epsilon$ to $\{10^{-4},10^{-4}\}$ respectively for FCG and $\{10^{-2},10^{-4}\}$ for FGMRES respectively. $2000$ data points were generated randomly in a unit cube for testing the convergence. The computational performance and convergence is shown in Fig. \ref{fig:fgmres}. The number of iterations of the preconditioned approaches are always less than those for the unpreconditioned cases. The computational cost per iteration is the least for CG compared to GMRES, FCG, and FGMRES. Incomplete LU (ILU) based preconditioners are marginally better in convergence (iterations) compared to our approach for better conditioned cases. But ILU (and other similar preconditioners) require explicit kernel matrix construction and rely on sparsity of the matrix to be solved and the absence of these properties in kernel matrices result in significantly higher computational cost compared to our preconditioners as well as the unpreconditioned solver. This makes them impractical.

We see that FCG needs increased accuracy of the inner linear system solution. In contrast, FGMRES is more forgiving of inner linear system error and only requires coarse accuracy to reduce the number of outer iterations to the same magnitude as FCG. On the other hand, especially for the ill-conditioned matrices, solving the inner Krylov method with fine accuracy takes much more time. Hence, given the ill-conditioned kernel matrices, the best FMGRES has the lesser number of outer iterations as well as smaller computation time.

The unpreconditioned algorithm of choice is CG, because of its lower storage and efficiency. However, FGMRES is the method of choice for preconditioned iterations. While GMRES requires extra storage in comparison to CG, FCG also requires this extra storage (for reorthogonalization), and we do not pay a storage penalty for our choice of FGMRES over FCG. In the sequel, we accordingly use FGMRES.

\section{Experiments\label{sec:Experiments}}
The preconditioner performance is illustrated on various datasets on different variants of kernel regression. We first look at GPR with a Gaussian kernel and then extend the preconditioned approach to a generalized (non-Gaussian) kernel. We also experiment on kriging \cite{AppliedGeostatistics} and report results on a large geostatistical dataset.

Although dataset-specific tuning of the preconditioner parameters can yield better results, this is impractical. We therefore use the following rules to set the preconditioner parameters.
\begin{itemize}
    \item The tolerance ($\epsilon$) for the preconditioner system solution is set at an order of magnitude larger than the outer iteration tolerance (e.g., outer tolerance $= 10^{-4}$, inner tolerance $=10^{-3}$).
    \item Similarly, the preconditioner regularizer $\delta$ is also set to an order of magnitude higher than the kernel regularizer $\gamma$. When the outer regularizer is $0$, the inner regularizer is set to $10^{-3}$
\end{itemize}
While this might not yield the best preconditioner system, it performs well in most cases from our experiments. In all experiments, the outer iteration tolerance was set to $10^{-6}$.

\subsection{Gaussian process regression (GPR): } The key computational bottleneck in GPR involves the solution of a linear system of kernel matrix. Direct solution via decompositions such Cholesky \cite{GPML_Rasmussen} requires O($N^2$) space and O($N^3$) time requirements. Alternatively, Mackay et al. \cite{GPML_Mackay} use \emph{CG} to solve the GPR in Eq. \ref{eq:GPR_mean}.

In our experiments, the covariance kernel parameters are estimated via maximum likelihood as in \cite{GPML_Rasmussen} with a small subset of the input data. We compare the performance of our preconditioner against a direct solution using \cite{GPML_Rasmussen}, our implementation of the CG approach in \cite{GPML_Mackay} and Incomplete LU based preconditioner on various standard datasets\footnote{\url{www.liaad.up.pt/~ltorgo/Regression/}}. The kernel matrix vector product in all compared scenarios was also accelerated using GPUML. Table \ref{table:fgmres_gpr} shows the corresponding result.

The convergence of the preconditioned FGMRES (both with ILU and our preconditioner) is consistently better than the unpreconditioned approach. Although for smaller datasets there is very little separating the computational performance of the solvers, the performance of our FGMRES with our preconditioner gets better for larger data sizes. This is because, for larger problems, cost per iteration in both CG and FGMRES increases, and thus a FGMRES which converges faster becomes significantly better than the CG-based approach. Further, for larger problems, both the direct method and ILU-preconditioning run into space issues due to the requirement of the physical construction of the kernel matrix.

\begin{table*}[bth]
\caption{\emph{Performance of FGMRES based Gaussian process regression against the direct, CG \cite{GPML_Mackay} and ILU-preconditioned solvers; $d$ is the dimension and $N$ is the size of the regression dataset with the Gaussian kernel. Total time taken for prediction is shown here, with the number of iterations for convergence indicated within parenthesis. The mean error in prediction between the two approaches was less than $10^{-6}$ in all the cases.}\label{table:fgmres_gpr}}
\centering
\begin{tabular}{||l|c|c|c|c||}
\hline\hline
Datasets ($d\times N$)	& Direct \cite{GPML_Rasmussen} & CG \cite{GPML_Mackay} 	& ILU & FGMRES\\
\hline\hline
\emph{Diabetes} ($3 \times 43$)	& $0.03$ 	& $\mathbf{0.03}$ ($8$) 	& $0.25$ ($4$) 	& $0.04$ (${3}$)\\
\emph{Boston Housing} ($14 \times 506$)	& $0.86$ 	& $0.67$ ($33$) 	& $0.86$ ($3$) 	& $\mathbf{0.62}$ (${3}$)\\
\emph{Pumadyn} ($9 \times 4499$)	& $63.61$ 	& $5.61$ ($32$) 	& $73.45$ ($4$) 	& $\mathbf{3.61}$ (${3}$)\\
\emph{Bank ($1$)} ($9 \times 4499$)	& $64.18$ 	& $6.53$ ($35$) 	& $74.73$ ($4$) 	& $\mathbf{4.28}$ (${3}$)\\
\emph{Robot Arm} ($9 \times 8192$)	& $232.61$ 	& $23.81$ ($75$) 	& $268.37$ ($3$) 	& $\mathbf{11.79}$ (${4}$)\\
\emph{Bank ($2$)} ($33 \times 4500$)	& $66.85$ 	& $49.40$ ($38$) 	& $76.54$ ($3$) 	& $\mathbf{37.74}$ (${3}$)\\
\emph{Census ($1$)} ($9 \times 22784$)	& $--$ 	& $117.45$ ($42$) 	& $--$ 	& $\mathbf{90.31}$ (${4}$)\\
\emph{Ailerons} ($41 \times 7154$)	& $170.76$ 	& $131.34$ ($31$) 	& $208.87$ ($2$) 	& $\mathbf{128.22}$ (${4}$)\\
\emph{2D Planes} ($11 \times 40768$)	& $--$ 	& $469.41$ ($31$) 	& $--$ 	& $\mathbf{415.30}$ (${6}$)\\
\emph{Census ($2$)} ($17 \times 22784$)	& $--$ 	& $663.70$ ($83$) 	& $--$ 	& $\mathbf{482.50}$ (${5}$)\\
\emph{Sarcos} ($28 \times 44484$)	& $--$ 	& $--$ 	& $--$ 	& $\mathbf{1090.85}$ (${4}$)\\
\hline\hline
\emph{\textbf{Kriging [Pacific Coast Data]}} & $--$ & $2,301\pm800$s & $--$ & $\mathbf{725\pm190}$s \\
($3 \times 179,065\pm35,405$) 	& $--$ & ($46\pm12$) & $--$ & ($3\pm1$)\\ \hline\hline
\end{tabular}
\end{table*}

Low rank approaches \cite{SnelsonGhahramani,Seeger2003,SPGP} also address the time complexity in kernel regression by working on an ``active set'' of set $M$ and reducing the time to O($M^2N$). We compared with the low rank GPR based on \cite{SPGP}, and found our approach to be superior. Because these approaches involve the solution of a large optimization problem, straightforward algorithmic acceleration or parallelization is not possible. Since the methods and accelerations used in this paper are significantly different from those in \cite{SPGP}, we have not reported these here.

To illustrate the applicability of our preconditioner to non-Gaussian kernels, we tested it on the Matern kernel\cite{GPR_Rasmussen},
\begin{equation}k(x_i,x_j)=(1 + \sqrt{3}d_{ij}) \exp(-\sqrt{3}d_{ij}),\label{eq:MaternKernel}\end{equation}
where $d_{ij}=\sqrt{\frac{\|x_i-x_j\|^2}{h^2}}$. We used the GPR framework for a binary classification problem and tested it on several standard datasets\footnote{\url{www.csie.ntu.edu.tw/~cjlin/libsvmtools/}}. The results are tabulated in Table \ref{table:fgmres_rlsc}. Here again, the FGMRES has a better computational performance than the other compared methods, thus illustrating its validity on non-Gaussian kernels.

\begin{table*}[bth]
\caption{\emph{Performance of FGMRES based Gaussian process regression against the direct, CG \cite{GPML_Mackay} and ILU-preconditioned solvers; $d$ is the dimension and $N$ is the size of the regression dataset with a non-Gaussian kernel (Matern). Total time taken for prediction is shown here, with the number of iterations for convergence indicated within parenthesis. The mean error in prediction between the two approaches was less than $10^{-6}$ in all the cases.}\label{table:fgmres_rlsc}}
\centering
\begin{tabular}{||l|c|c|c|c||}
\hline\hline
Datasets ($d\times N$)	& Direct \cite{GPML_Rasmussen} & CG \cite{GPML_Mackay} 	& ILU & FGMRES\\
\hline\hline
\emph{Heart} ($14 \times 270$)	& $0.24$ 	& $0.10$ ($7$) 	& $0.31$ ($2$) 	& $\mathbf{0.09}$ (${2}$)\\
\emph{Iris} ($5 \times 150$)	& $\mathbf{0.08}$ 	& $0.10$ ($28$) 	& $0.10$ ($3$) 	& $0.10$ (${3}$)\\
\emph{Sonar} ($61 \times 208$)	& $0.17$ 	& $\mathbf{0.06}$ ($2$) 	& $0.46$ (${3}$) 	& $0.16$ (${2}$)\\
\emph{Diabetes} ($9 \times 768$)	& $1.83$ 	& $0.50$ ($23$) 	& $2.06$ (${2}$) 	& $\mathbf{0.35}$ ($3$)\\
\emph{Glass} ($10 \times 214$)	& $\mathbf{0.15}$ 	& $0.32$ ($45$) 	& $0.21$ ($3$) 	& $0.36$ (${4}$)\\
\emph{German} ($25 \times 1000$)	& $3.26$ 	& $\mathbf{0.33}$ ($4$) 	& $3.92$ ($3$) 	& $0.45$ (${2}$)\\
\emph{Australian} ($15 \times 690$)	& $1.59$ 	& $0.71$ ($25$) 	& $1.95$ ($3$) 	& $\mathbf{0.53}$ (${3}$)\\
\emph{Vehicle} ($19 \times 846$)	& $2.80$ 	& $0.66$ ($16$) 	& $2.76$ ($4$) 	& $\mathbf{0.55}$ (${3}$)\\
\emph{Splice} ($61 \times 1000$)	& $4.51$ 	& $\mathbf{0.39}$ ($25$) 	& $5.01$ ($3$) 	& $0.93$ (${2}$)\\
\emph{Fourclass} ($3 \times 862$)	& $3.26$ 	& $2.56$ ($204$) 	& $2.45$ ($2$) 	& $\mathbf{1.88}$ (${3}$)\\
\emph{letter} ($17 \times 15000$)	& $119.51$ 	& $186.58$ ($35$) 	& $129.51$ ($3$) 	& $\mathbf{69.54}$ (${4}$)\\
\hline\hline
\end{tabular}
\end{table*}

\subsection{Kriging: }
\begin{figure}[bth]
\centering
\includegraphics[width=\linewidth]{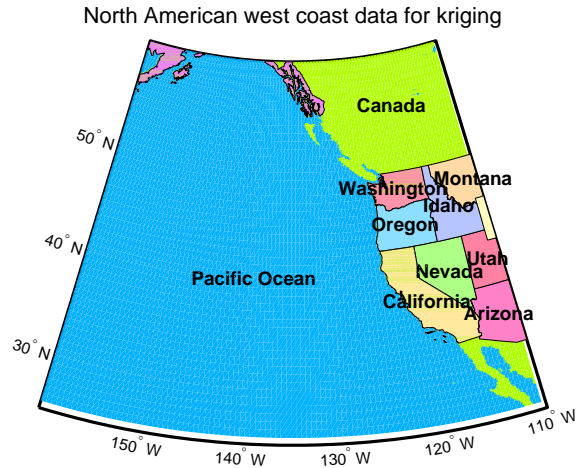}
\caption{\emph{Kriging was performed on the data recorded along the west coast of North America shown here}\label{fig:KrigingMap}}
\end{figure}

We compared FGMRES-based kriging against the CG version on the ocean chlorophyll concentration data recorded along the Pacific coast of North America (the data map is shown in Fig. \ref{fig:KrigingMap}) obtained from National Oceanic and Atmospheric Administration\footnote{\url{http://coastwatch.pfel.noaa.gov/}}. We look at the $7$-day aggregate of the chlorophyll concentration, which is recorded on a grid of $416\times600$. However, this includes several locations with missing data or those located over land. This results in approximately $179,065\pm35,405$ data samples per week.

It was observed that the CG-based approach converges in $46\pm12$ iterations in $2,301\pm800$s, whereas, the FGMRES converges in just $\mathbf{3\pm1}$ (outer) iterations in $\mathbf{725\pm190}$s, resulting in over $3$X speedup.

\section{Conclusions and discussions\label{sec:conclusions}}
A method to improve convergence of Krylov methods used in kernel methods was demonstrated. The key contributions of the paper are as follows,
\begin{itemize}
\item A novel yet simple preconditioner is proposed to solve a linear system with a kernel matrix using flexible Krylov methods.
\item A technique to accelerate the inner preconditioner system using truncated CG with fast matrix vector products was developed.
\item Rules to select the preconditioner parameter were shown.
\end{itemize}
The core preconditioning strategy proposed here will soon be released as an open source package with Matlab bindings.


\end{document}